\documentclass{article}
\usepackage{amssymb,latexsym,amsmath,amsthm,mathrsfs}
\usepackage{graphicx}
\newcommand{\comment}[1]{}
\begin{document}
\title{A more accurate treatment of the problem of drawing the shortest line on a surface\footnote{Presented to the St. Petersburg Academy on
January 25, 1779.
Originally published as
{\em Accuratior evolutio problematis de linea brevissima in superficie quacunque
ducenda},
Nova acta academiae scientiarum Petropolitanae \textbf{15} (1806),
44--54.
E727 in the Enestr{\"o}m index.
Translated from the Latin by Jordan Bell,
Department of Mathematics, University of Toronto, Toronto, Canada.
Email: jordan.bell@utoronto.ca}}
\author{Leonhard Euler}
\date{}
\maketitle

1. For a surface on which a shortest line is to be drawn, given
this differential equation between the three orthogonal coordinates
$x,y,z$: $dz=fdx+gdy$, where
$f$ and $g$ are functions of both $x$ and $y$, let it be
\[
df=\alpha dx+\beta dy \quad \textrm{and} \quad dg=\beta dx+\gamma dy.
\]
Having done this, since an element of any line drawn on this surface would be
\[
\sqrt{dx^2+dy^2+dz^2},
\]
with the previous value put in place of $dz$ an element of this curve will be
\[
=\sqrt{dx^2+dy^2+(fdx+gdy)^2};
\]
whence if we put $dy=pdx$, this element will be $dx\sqrt{1+pp+(f+gp)^2}$.

2. This integral formula, which is to be brought to a minimum, will be
\[
\int dx\sqrt{1+pp+(f+gp)^2},
\]
which
I indicated in general by $\int Zdx$
in my Treatise
{\em Methodus inveniendi lineas curvas Maximi Minimive proprietate gaudentes},
so that in this case it would be
\[
Z=\sqrt{1+pp+(f+gp)^2}.
\]
Then, having put $dZ=Mdx+Ndy+Pdp$,
I have shown that the nature of a Minimum or a Maximum is such that it is 
expressed by this equation:\footnote{Translator: This is namely
the Euler-Lagrange equation in the calculus of variations.} $Ndx=dP$, which clearly leads to differentials of
the second order.

3. Therefore since
$Z^2=1+pp+(f+gp)^2$,
let us differentiate this formula and let the elements be separated into
the three types,
namely
$dx,dy,dp$, and 
then be represented in this way:
\[
ZdZ=dx(\alpha+\beta p)(f+gp)+dy(\beta+\gamma p)(f+gp)+dp(p+g(f+gp)).
\]
Then since I have put in general $dZ=Mdx+Ndy+Pdp$, we will have in
this case:
\begin{eqnarray*}
M&=&\frac{(\alpha+\beta p)(f+gp)}{Z}\\
N&=&\frac{(\beta+\gamma p)(f+gp)}{Z}\\
P&=&\frac{p+g(f+gp)}{Z}.
\end{eqnarray*}
Hence (since $\beta dx+\gamma pdx=dg$) it becomes
\[
Ndx=\frac{dg(f+gp)}{Z},
\]
whence the equation for our sought curve will be
\[
\frac{dg(f+gp)}{Z}=d\cdot\frac{p+g(f+gp)}{Z}.
\]
For expanding this equation, for the sake of brevity let us put
$p+g(f+gp)=S$, and we will have:
\[
\frac{dg(f+gp)}{Z}=\frac{dS}{Z}-\frac{SdZ}{ZZ}
\]
or
\[
dg(f+gp)=dS-\frac{SdZ}{Z}.
\]
Therefore because $dS=dp+dg(f+gp)+gd(f+gp)$, our equation will be
\[
0=dp+gd(f+gp)-\frac{SdZ}{Z}.
\]
But on the other hand it is:
\[
\frac{dZ}{Z}=\frac{pdp+(f+gp)d(f+gp)}{1+pp+(f+gp)^2},
\]
which should be multiplied by $S=p+g(f+gp)$. Then multiplying by the denominator
$1+pp+(f+gp)^2$ we will have:
\[
0=dp+(g-fp)d(f+gp)-gpdp(f+gp)+dp(f+gp)^2
\]
or
\[
0=dp+(g-fp)d(f+gp)+fdp(f+gp),
\]
which we can then write in this form:
\[
0=dp(1+ff+gg)+(g-fp)(df+pdg).
\]

4. Although this equation is simple enough, it is not clear however how
it could be reduced to differentials of the first degree. But I have
observed that the problem can be dealt with by the following substitution,
namely: $v=\frac{g-fp}{f+gp}$; whence it will be $p=\frac{g-fv}{gv+f}$,
and now by differentiating we deduce
\[
dp=\frac{-(ff+gg)dv+(1+vv)(fdg-gdf)}{(f+gv)^2}.
\]
Further, it will be 
\[
g-fp=\frac{v(ff+gg)}{f+gv},
\]
and next
\[
df+pdg=\frac{fdf+gdg+v(gdf-fdg)}{f+gv}.
\]
With these substituted in, the equation becomes:
\begin{eqnarray*}
0&=&-dv(ff+gg)(1+ff+gg)+v(ff+gg)(fdf+gdg)\\
&&+(1+vv)(fdg-gdf)+(ff+gg)(fdg-gdf).
\end{eqnarray*}

5. In order to turn this into a simpler equation let us set $ff+gg=hh$,
and it will be $fdf+gdg=hdh$, and next let $\frac{g}{f}=k$, so that it
would become $fdg-gdf=ffdk$, and thus our equation
is contracted into this form:
\[
0=-hhdv(1+hh)+h^3vdh+(1+hh+vv)ffdk.
\]
Also since $g=fk$, it will be $ff(1+kk)=hh$ and then $ff=\frac{hh}{1+kk}$,
from which we will have:
\[
0=-dv(1+hh)+vhdh+(1+hh+vv)\frac{dk}{1+kk};
\]
next by putting $v=s\sqrt{1+hh}$ this equation is reduced to this form:
\[
0=-ds\surd(1+hh)+\frac{dk(1+ss)}{1+kk}.
\]
Therefore now the quantity $s$ can be exhibited separately from the others,
as
\[
\frac{ds}{1+ss}=\frac{dk}{(1+kk)\sqrt{1+hh}},
\]
which seems to be the simplest form to which it can be brought in general.

6. While we have attempted to determine the two variables $y$ and $x$
by the single variable $z$, since all three are involved in equal measure
in the calculation, it appeared to me to treat this whole question so that
all the formulae involve in equal measure the three coordinates $x,y,z$;
this idea will be put to good use in the following investigations that I will take up.

\begin{center}
{\Large Supplement}
\end{center}

7. Let the differential equation of a given surface be:
$pdx+qdy+rdz=0$,
where $p,q,r$ are functions of the coordinates $x,y,z$:
whence for this equation to be possible, this condition must be
satisfied:\footnote{Translator: cf. p. 101 of D. J. Struik,
{\em Outline of a history of differential geometry: I},
Isis \textbf{19} (1933), no. 1, 92--120}
\[
\frac{pdq-qdp}{dz}+\frac{qdr-rdq}{dx}+\frac{rdp-pdr}{dy}=0.
\]
With this done, the following equation for drawing the shortest line
on this surface will be obtained, which involves the three coordinates
$x,y,z$ is equal measure:
\[
ddx(qdz-rdy)+ddy(rdx-pdz)+ddz(pdy-qdx)=0.
\]
Or if we put for the sake of brevity:
\begin{eqnarray*}
dyddz-dzddy&=&f,\\
dzddx-dxddz&=&g,\\
dxddy-dyddx&=&h,
\end{eqnarray*}
it will be $fp+gq+hr=0$, then indeed on the other hand
$fdx+gdy+hdz=0$.
Next if an element of the shortest curve is put $=ds$, it will be
$ds^2=dx^2+dy^2+dz^2$, then indeed also
\[
\frac{dds}{ds}=\frac{qddz-rddy}{qdz-rdy}=\frac{rddx-pddz}{rdx-pdz}=
\frac{pddy-qddx}{pdy-qdx}.
\]

\begin{center}
{\Large Application to a spherical surface}
\end{center}

8. Let the equation for the surface be $xdx+ydy+zdz=0$, so that we have
here $p=x,q=y,r=z$, and the first equation for the shortest path
will be the following:
\[
ddx(ydz-zdy)+ddy(zdx-xdz)+ddz(xdy-ydx)=0,
\]
whose complete integral is thus $\alpha x+\beta y+\gamma z=0$, which is
apparent from the nature of the matter. The question is therefore
reduced to how this integral can be worked out.

9. Now with the earlier equation $fx+gy+hz=0$, if for this equation
we put $\Pi=\frac{zdx-xdz}{ydx-xdy}$, it will be
$d\Pi=d\cdot \frac{zdx-xdz}{ydx-xdy}$ and then
\[
d\Pi=\frac{zddx-xddz}{ydx-xdy}-\frac{(zdx-xdz)(yddx-xddy)}{(ydx-xdy)^2},
\]
or by expanding
\[
d\Pi=\frac{x}{(ydx-xdy)^2}((dyddz-dzddy)x+(dzddx-dxddz)y+(dxddy-dyddx)z)
\]
and with the $f,g,h$ that have been introduced it will be
\[
d\Pi=x\frac{(fx+gy+hz)}{(ydx-xdy)^2}.
\]
Since on the other hand it is $fx+gy+hz=0$, it will be $d\Pi=0$ and hence
$\Pi$ is a constant quantity, which if we put $=A$, the differential
equation of the first degree $\Pi=\frac{zdx-xdz}{ydx-xdy}$ can be
thus expressed:
\[
A(ydx-xdy)=zdx-xdz,
\]
which when divided by $xx$ will be integrable; for it would become
\[
\frac{Ay}{x}=\frac{z}{x}+B \quad \textrm{or} \quad Ay-Bx-z=0
\]
or with the constants switched
\[
\alpha x+\beta y+\gamma z=0.
\]
Since this equation is clearly everywhere drawn from the center
of the sphere, great circles arise on the surface of the sphere;
whence it follows that all great circles are all the shortest paths which can
be drawn on the surface of the sphere.

10. Since in these calculations everything is typically reduced to a single
variable, if for effecting this we put $dy=tdx$ and $dx=udx$,
taking $dx$ as constant, the first equation will be as follows:
\[
dt(r-pu)+du(pt-q)=0.
\]
And the equation for the surface will be $p+qt+ru=0$;
whence, since it then becomes $p=-qt-ru$, the former equation takes this form:
\[
dt(r+qtu+ruu)-du(q+rtu+qtt)=0.
\]
Next it will be
\[
f=dx^2(tdu-udt), \quad g=-dx^2du, \quad h=dx^2dt,
\]
then indeed $ds^2=dx^2(1+tt+uu)$ and finally
\[
\frac{dds}{ds}=\frac{tdt+udu}{1+tt+uu}=\frac{qdu-rdt}{qu-rt}=
-\frac{pdu}{r-pu}=\frac{pdt}{pt-q}.
\]

11. And if desired to introduce as a fourth variable the angle $\varphi$,
by putting $dx=td\varphi,dy=ud\varphi,dz=vd\varphi$, the equation
for the surface will be
\[
pt+qu+rv=0.
\]
Next for the letters $f,g,h$ we will have
\begin{eqnarray*}
f&=&d\varphi^2(udv-vdu),\\
g&=&d\varphi^2(vdt-tdv),\\
h&=&d\varphi^2(tdu-udt),
\end{eqnarray*}
thus it will then be $ft+gu+hv=0$. The equation for the shortest path will be:
\[
fp+gq+hr=p(udv-vdu)+q(vdt-tdv)+r(tdu-udt)=0,
\]
and it would finally become $ds^2=d\varphi^2(tt+uu+vv)$ and hence
\[
\frac{dds}{ds}=\frac{tdt+udu+vdv}{tt+uu+vv}=\frac{qdv-rdu}{qv-ru}=
\frac{rdt-pdv}{rt-pv}=\frac{pdu-qdt}{pu-qt}.
\]

12. Since no way presents itself to us for integrating the general equation
given above for 
drawing the shortest path on a surface, even though many cases can be
given in which integration of the equation for the curve succeeds, 
it will be worthwhile to have expanded several of these here as a conclusion.

13. Let us begin with the case in which one of the quantities $p,q,r$ vanishes.
In particular, if $r=0$, then the equation for the surface
will be $pdx+qdy=0$, in which case therefore the surface is a cylinder, whose
base is determined by the equation $pdx+qdy$. And putting
in the equation for $\frac{dds}{ds}$ the given $r=0$, it would become
$\frac{dds}{ds}=\frac{ddz}{dz}$, whose integral is
$lds=ldz+l\alpha$ and hence, taking numbers, $ds=\alpha dz$ or
\[
dx^2+dy^2+dz^2=\alpha \alpha dz^2 \quad \textrm{or} \quad
dx^2+dy^2=dz^2(\alpha\alpha-1).
\]
It will then be
\[
z\sqrt{\alpha\alpha-1}=\int\sqrt{dx^2+dy^2},
\]
where $\int \sqrt{dx^2+dy^2}$ expresses an element
of the base of the curve; whence it is apparent that the altitude
$z$ of a cylinder is always proportional to the arc of the base.

14. Let us consider the case in which $p=x$ and $q=y$, where the
equation for the surface will thus be:
\[
xdx+ydx+rdz=0,
\]
in which all spherical bodies or rotated surfaces are contained.
Then it will further be
\[
\frac{dds}{ds}=\frac{xddy-yddx}{xdy-ydx},
\]
whose integral is
\[
lds=l(xdy-ydx)+la,
\]
and hence by taking numbers
\[
\frac{ds}{a}=xdy-ydx \quad \textrm{and so} \quad \frac{dx^2+dy^2+dz^2}{aa}=(xdy-ydx)^2
\]
or switching the constants
\[
(dx^2+dy^2+dz^2)AA=(xdy-ydx)^2.
\]
For integrating this equation again let us put $x=v\cos\varphi$ and
$y=v\sin\varphi$, and it will be
\[
dx^2+dy^2=dv^2+vvd\varphi^2;
\]
then indeed for the curve it will be $vdv+rdz=0$, where $r$ is some function
of $v$, so that $dz=-\frac{vdv}{r}$. Next indeed it will be
\[
xdy-ydx=vvd\varphi;
\]
and with all these substitutions our equation will be:
\[
AA(dv^2+vvd\varphi^2+\frac{vvdv^2}{rr})=AAds^2=v^4d\varphi^2,
\]
from which in turn we gather
\[
d\varphi^2=\frac{AA(dv^2+\frac{vvdv^2}{rr})}{v^4-vvAA}=\frac{AAdv^2(rr+vv)}{rrvv(vv-AA)}
\]
and hence
\[
d\varphi=\frac{Adv}{rv}\sqrt{\frac{rr+vv}{vv-AA}}.
\]

15. For other cases our general equation above can be put to more use.
And first indeed, because so far things have depended on a ratio between
the quantities $p,q,r$, one of
which will be assumed at our pleasure.
Thus let $r=-1$, so it becomes $dz=pdx+qdy$, and let $p$ and $q$ be functions
of $x$ and $y$ with it being $(\frac{dp}{dy})=(\frac{dq}{dz})$. Next
let us put $dy=\pi dx$, and it will be $dz=(p+\pi q)dx$.
Then, taking the element $dx$ as constant, so that it would be $ddx=0$,
it will be
\[
ddy=d\pi dx \quad \textrm{and} \quad ddz=(dp+\pi dq+qd\pi)dx.
\]
Now the three letters $f,g,h$ can be expressed in the following
way:
\begin{eqnarray*}
f&=&dx^2(\pi dp-pd\pi +\pi\pi dq),\\
g&=&-dx^2(dp+\pi dq+qd\pi),\\
h&=&d\pi dx^2.
\end{eqnarray*}
Indeed we have seen that the equation for the shortest path is
\[
pf+gq+hr=0,
\]
which therefore takes this form:
\[
-d\pi(1+pp+qq)+dp(\pi p-q)+\pi dq(\pi p-q)=0
\]
or
\[
d\pi(1+pp+qq)+(dp+\pi dq)(q-\pi p)=0.
\]

16. Since the two formulae $p+\pi q$ adn $q-\pi p$ are the central
components of this equation, it will be very helpful to work out
a relation between them. 
To this end let us put $\frac{q-\pi p}{p+\pi q}=v$, whence it is now
$\pi=\frac{q-vp}{p+vq}$; then indeed in turn
\[
g-\pi p=\frac{v(pp+qq)}{p+vq},
\]
and on the other hand it will be
\[
dp+\pi dq=\frac{pdp+qdq+v(qdp-pdq)}{p+qv}.
\]
Now if we put $q=up$, it will be
\[
\pi=\frac{u-v}{1+uv} \quad \textrm{and then} \quad d\pi=\frac{du(1+vv)-dv(1+uu)}{(1+uv)^2}.
\]
Next let us put $pp+qq=tt$, and since we let $q=up$, it will be
\[
pp=\frac{t}{1+uu}\quad \textrm{and} \quad d\cdot\frac{p}{q}=du=\frac{pdq-qdp}{pp},
\]
and hence
\[
pdq-qdp=ppdu=\frac{ttdu}{1+uu},
\]
and when these values are substituted, because
\[
q-\pi p=\frac{vtt}{p(1+vu)} \quad \textrm{and} \quad dp+\pi dq=
\frac{tdt-(vttdu):(1+uu)}{p(1+uv)},
\]
it will be
\[
0=-\frac{(du(1+vv)-dv(1+uu))(1+tt)}{(1+uv)^2}-\frac{vtt(t(1+uu)dt-vttdu)}{pp(1+uu)(1+uv)^2}
\]
or
\[
(1+tt)(du(1+vv)-dv(1+uu))+vt((1+uu)dt-vtdu)=0,
\]
which can then be reduced to this form:
\[
du((1+vv)(1+tt)-vvtt)-dv(1+tt)(1+uu)+vtdt(1+uu)=0
\]
or this even more neat one:
\[
\frac{du}{1+uu}(1+vv+tt)-dv(1+tt)+vtdt=0.
\]
Now let us put $v=w\sqrt{1+tt}$, and it will be
\[
dw=\frac{dv(1+tt)-vtdt}{(1+tt)^{\frac{3}{2}}}
\]
that is, it will be
\[
dv(1+tt)-vtdt=(1+tt)^{\frac{3}{2}}dw;
\]
then indeed it will be
\[
1+tt+vv=(1+tt)(1+ww),
\]
and with these values substituted our equation itself thus becomes
\[
\frac{du}{1+uu}(1+tt)(1+ww)-(1+tt)^{\frac{3}{2}}dw=0,
\]
hence by separating we obtain
\[
\frac{du}{1+uu}=\frac{dw\sqrt{1+tt}}{1+ww},
\]
and consequently
\[
\frac{dw}{1+ww}=\frac{du}{(1+uu)\sqrt{1+tt}},
\]
which equation can always be intergated whenever $t$ is a function
of $u$ or whenever $pp+qq$ is a function of $\frac{q}{p}$
or $q$ a function of $p$.

17. It will further turn out, so that $q$ is a function of $p$, first
if $z$ and $y$ are thus determined by $x$ and another new variable $\omega$,
so that it would be $y=Ax$ and $z=Bx$, with $A$ and $B$ being some functions
of $\omega$.
Therefore since we have put $dz=pdx+qdy$, it will be
\[
Bdx+xdB=pdx+qAdx+qxdA,
\]
where the terms involving the differential should be compared
separately, whence it will be $p=B-Aq$;
and comparing separately the terms containing the quantity $x$ it will
be $q=\frac{dB}{dA}$ and hence $p=\frac{BdA-AdB}{dA}$
And thus
$p$ and $q$ are functions of $\omega$ and hence $tt=pp+qq$ and $u=\frac{p}{q}$
will be functions of the single quantity $\omega$ and $\sqrt{1+tt}$ will be a function of $u$.
Because of this, the equation found above for the shortest path admits
integration.
Moreover, in this case, namely in which $y=Ax$ and $z=Bx$, 
a conical surface follows constructed upon any base.

18. The equation given above would then be integrable when we put
$y=Ax+C$ and $z=Bx+D$; for then it will be
\[
dz=pdx+qdy=Bdx+xdB+dD
\]
and because $dy=Adx+xdA+dC$, it will also be
\[
dz=pdx+qdy=pdx+Aqdx+xqdA+qdC
\]
and hence, comparing the members containing the quantity $x$ with each other,
and then indeed those which are endowed with the differential $dx$, it will be
\[
B=p+Aq \quad \textrm{and} \quad dB=qdA,
\]
then
\[
q=\frac{dB}{dA} \quad \textrm{and} \quad p=\frac{BdA-AdB}{dA}.
\]
Besides these indeed it should be $dD=qdC=\frac{dBdC}{dA}$,
or the functions $A,B,C,D$ should thus be compared that
$dAdD=dBdC$. If this occurs, $p$ and $q$ will again be function of the
single variable $\omega$ and hence $\sqrt{1+tt}$ will also be a function
of $u$, in which case too it is possible to define a shortest path.
This case seems indeed to complete all the planar surfaces which can be
explained in the plane.

\end{document}